\documentstyle[12pt,amsfonts,epsf]{article}

\newcommand{\cliff}{{\cal C}}

\newtheorem{lemma}{Lemma}
\newtheorem{exercise}[lemma]{Exercise}

\newtheorem{theorem}[lemma]{Theorem}
\newtheorem{corollary}[lemma]{Corollary}
\newtheorem{proposition}[lemma]{Proposition}
\newtheorem{definition}[lemma]{Definition}

\begin{document}

\title{Rotations of the three-sphere and symmetry of the Clifford Torus}
\author{John McCuan\thanks{Work supported in part by the University of California, Berkeley.} \ \ and\ \  Lafe Spietz}


\maketitle

\begin{abstract}
We describe decomposition formulas for rotations of ${\Bbb R}^3$ and 
${\Bbb R}^4$ that have special 
properties with respect to stereographic projection.  We use the lower 
dimensional decomposition to analyze stereographic projections of great 
circles in ${\Bbb S}^2 \subset {\Bbb R}^3$.  This analysis provides a 
pattern for our analysis of stereographic projections of the Clifford 
torus ${\cal C}\subset {\Bbb S}^3 \subset{\Bbb R}^4$.  We use the higher 
dimensional decomposition to prove a symmetry assertion for stereographic 
projections of ${\cal C}$ which we believe we are the first to observe 
and which can be used to characterize the Clifford torus among embedded 
minimal tori in ${\Bbb S}^3$---though this last assertion goes beyond 
the scope of this paper.  
An effort is made to intuitively motivate all necessary 
concepts including rotation, stereographic projection, and symmetry.
\end{abstract}

\pagestyle{headings}
{\makeatletter
\gdef\@evenhead{\thepage\hfil\small\sc John McCuan and Lafe Spietz}%
\gdef\@oddhead{%
\ifodd\thepage
\hfil{\small\uppercase{Three-sphere rotations and the Clifford Torus}}\hfil
\thepage
\else{\thepage\hfil\small\sc \uppercase{John McCuan and Lafe Spietz}\hfil}\fi}
}


\section*{Introduction}

It is known (and intuitively believable) that the spheres in ${\Bbb R}^3$ 
are characterized by being the only compact surfaces with planes of 
symmetry whose normals exhaust all possible directions.  That is, 
{\em if $S$ is a compact surface and, for each unit vector ${\bold n}$ in 
${\Bbb R}^3$, there is some plane $\Pi$ with normal ${\bold n}$ such 
that $S$ is invariant under reflection in $\Pi$, then $S$ is a sphere.}  
Technically, $S$ could be some collection of concentric spheres, but this 
can be fixed by requiring explicitly that $S$ be connected, i.e., a single 
surface.  The point is that spheres can be characterized by their 
(reflectional) symmetry.

In this paper we describe an analogous symmetry condition for a certain 
surface, the Clifford torus.  Our task is complicated by 
the fact that this torus is located in the three-sphere, ${\Bbb S}^3$, 
where our intuition from ${\Bbb R}^3$ is of limited use.  For this 
reason we will employ a certain transformation, {\em stereographic 
projection,} that allows us to realize ${\Bbb S}^3$ (at least most of it) 
in the Euclidean space ${\Bbb R}^3$.  In fact, our symmetry condition 
will apply, more properly, to the stereographic projections of the 
Clifford torus.  Furthermore, our symmetry condition is substantially more 
complicated than the simple one above for spheres, and we discuss at some 
length why it is a natural one.

Another objective of the paper is to give what we consider novel, 
geometrically based, expositions of several well known topics.  (Some 
of these are mentioned briefly below.)  From this point of view, we offer 
an introduction to ${\Bbb S}^3$ that we hope is a geometric counterpart 
to the algebraic treatment in, for example, \cite{ZulCha}.

The paper is organized as follows.  In the following \S\ref{SP} we 
review stereographic projection of ${\Bbb S}^2$ and discuss a 
decomposition formula for rotations of ${\Bbb S}^2 \subset {\Bbb R}^3$.  
We show, in particular, that projections of rotations of great circles 
are circles in ${\Bbb R}^2$ whose size and position are given in terms of 
the parameters of the decomposition.  This discussion is somewhat 
artificial because it is easy to show that the projection of essentially 
any circle in ${\Bbb S}^2$ is a circle in ${\Bbb R}^2$ whose center and 
radius are easy to calculate.  This is the case, however, owing to the 
fact that a circle in ${\Bbb S}^2$ is the intersection of an affine 
subspace (a plane) with ${\Bbb S}^2$.  This luxury is not afforded us 
by the Clifford torus, and the decomposition technique presented here 
will be used with considerable advantage in the more complicated 
higher dimensional case.  Furthermore, our discussion of stereographic 
projection of ${\Bbb S}^2$ is used in \S\ref{SYM} to give an exposition 
of symmetry for planar sets, and it provides intuition for the higher 
dimensional stereographic projection considered in \S\ref{MR}.

We also give 
considerable attention to 
building up intuition about {\em rotations}, especially rotations of 
${\Bbb R}^4$ since we consider ${\Bbb S}^3$ as a subset of ${\Bbb R}^4$.  
Recall that 
rotations of ${\Bbb R}^2$ (centered at zero) can be represented by 
matrices of the form
\[
	\left({\begin{array}{cc}
		\cos\theta & -\sin\theta \\
		\sin\theta & \ \ \cos\theta \end{array} } \right).
\]
That is, given a rotation $R: {\Bbb R}^2 \to {\Bbb R}^2$ there is a 
matrix $M$ of the form given above such that for each vector 
$\bold x$, we have that $R({\bold x})$ is given by the matrix multiplication 
$M{\bold x}$.  Through such a representation we are immediately presented 
with two algebraic facts, namely that $R$ is {\em linear} and that 
${\rm det}\, M = 1$.  It is also easy to check from this representation 
that $R$ is {\em orthogonal}, i.e., it preserves the orthonormality of 
bases.  Many authors define a 
{\em rotation} of ${\Bbb R}^3$ to be an orthogonal linear transformation 
corresponding 
to a matrix of determinant 1.  This definition is concise and 
computationally convenient, but we find it unintuitive.  It can perhaps 
be argued (via the parallelogram rule) that linearity is an intuitive 
assumption, but the role played by the determinant is difficult to see 
geometrically for rotations of ${\Bbb R}^3$, much less for rotations 
of  ${\Bbb R}^4$. 
In an appendix to this paper we take the point 
of view that rotations are {\em rigid motions\/} (i.e., distance 
preserving transformations) that fix the origin and result from 
{\em smooth ``homogeneous'' motions.\/}    
We then prove linearity, orthogonality, and 
representation by matrices of determinant 1 (thereby showing the two 
definitions are equivalent).  This appendix may be read at any time, but 
we recommend reading it after \S\ref{SP} and before \S\ref{MR}.

In \S\ref{SYM} we introduce {\em circles of Apollonius\/} and show that 
their symmetry as a family of planar curves is a natural generalization 
of the symmetry exhibited by 
concentric circles.  We incorporate in this discussion an explanation of 
why symmetry and reflection about circles are natural generalizations of 
symmetry and reflection about lines.  Furthermore, we observe that the 
generalized symmetry exhibited by circles of Apollonius is 
described naturally in terms of a {\em line of centers\/} $l$. 
More precisely, given a family of circles of Apollonius ${\cal A}$, we 
can find, for each point ${\bold x}$ on $l$, an {\em orthogonal 
Steiner circle\/} $C$ with center ${\bold x}$ and the property that every 
circle in ${\cal A}$ is symmetric with respect to $C$.  Finally, we 
observe that this formulation of generalized symmetry can be easily 
extended to surfaces in ${\Bbb R}^3$.  This sets the stage for our 
main (and we believe original) result which is roughly as follows.  
{\em Given any stereographic 
projection $Q$ of the Clifford torus, there is a line $l$ in 
${\Bbb R}^3$, and for each point ${\bold x}$ on $l$ there is 
a sphere $S$ centered at ${\bold x}$ so that $Q$ is symmetric with 
respect to $S$.\/}  See Figure~\ref{figcliff}.  A careful statement 
and proof, which include 
conditions on the radii of the spheres of symmetry, are given in \S\ref{MR}.  
The proof follows, in outline, the discussion of \S\ref{SP} and \S\ref{SYM}.

From a wider perspective, the Clifford torus and its stereographic 
projections are interesting surfaces primarily due to curvature 
considerations that are beyond the scope of this paper.  More primitively, 
they are interesting because they are critical points for 
certain functionals---which 
they are {\em believed\/} to minimize.  For further information, see 
\cite{BryDua,HsuMin,WilRie,YauSem}.  In this framework, the Clifford torus 
plays a role in ${\Bbb S}^3$ similar to that of the plane {\em and} 
the sphere in ${\Bbb R}^3$.  It has Gauss curvature and mean curvature 
zero like the plane, and it has constant mean curvature and is 
compact like the sphere.  The symmetry properties shown in this 
paper may be considered as a first step in getting a feel for 
curvature of surfaces in ${\Bbb S}^3$, and while the main result stated 
above may seem curious at first sight, it is completely analogous to 
the observation that any sphere in ${\Bbb R}^3$ has a plane of reflective 
symmetry with any given normal direction.

We extend our thanks to Ed Bueler, David Hoffman, and Silvio Levy who 
gave us useful comments that greatly improved the exposition.

\section{Stereographic Projection}\label{SP}
In this section we consider the two-dimensional sphere ${\Bbb S}^2$ as 
a subset of ${\Bbb R}^3$:  
\[
	{\Bbb S}^2 = \{{\bold x} = (x,y,z): |{\bold x}| = 1\}.
\]  
As usual $|{\bold x}| = \sqrt{x^2+y^2+z^2}$.

We are interested in a map called {\em stereographic projection} that sends 
${\Bbb S}^2$ (except for one point) into the $x,y$-plane, and we 
are interested 
in how the image of a certain set changes as the sphere is rotated.  To 
be precise, stereographic projection $\pi$ given by 
\begin{equation}\label{sp}
	\pi({\bold x}) = {1\over{1-z}}(x,y)
\end{equation}
maps ${\Bbb S}^2\backslash\{(0,0,1)\}$ onto ${\Bbb R}^2 = \{ (x,y,0) \}$ in a 
one-to-one fashion.

\begin{exercise}\label{proj.ex}
	Show that $\pi$ is one-to-one and onto.
\end{exercise}

\noindent The map (\ref{sp}) has a convenient geometric interpretation:  
\begin{quote}
	Each 
	point ${\bold p}\in {\Bbb S}^2\backslash\{(0,0,1)\}$ determines a 
	unique line passing 
	through ${\bold p}$ and $(0,0,1)$.  
	The line $l$, in turn, intersects the $x,y$-plane in a 
	unique point ${\bold q}$.  We set $\pi({\bold p}) = {\bold q}$.
\end{quote}
\begin{figure}[hb]
\centerline{
\begin{picture}(220,100)
	\put(76,82){${\bold p}$}
	\put(34,36){${\bold q}$}
	\put(85,80){\circle*{3}} 
	\put(25,40){\circle*{3}} 
	\put(109,96){\circle*{3}}
	\put(0,0){\epsfxsize=3in\epsffile{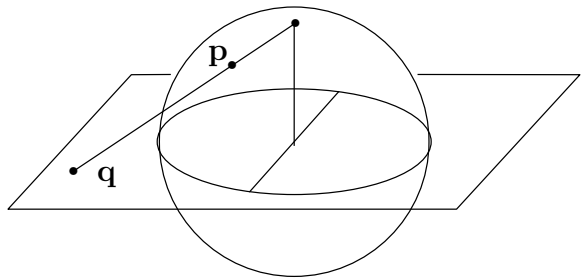}}
\end{picture}}
\centerline{ }
\caption{Stereographic Projection.}
\label{fig0}
\end{figure}
\begin{exercise}\label{ext.ex}
\begin{description}
\item[(i)] Using this geometric statement as a definition, derive formula 
	{\rm (\ref{sp})} for $\pi$.
\item[(ii)] The geometric definition makes sense for any point in ${\Bbb R}^3$ 
	with third component not equal to 1.  In this way $\pi$ may be 
	extended from ${\Bbb S}^2\backslash\{(0,0,1)\}$ to 
	${\Bbb R}^3\backslash\{(x,y,1)\}$.  We call this extended map 
	$\bar{\pi}$.  Does formula {\rm (\ref{sp})} still apply?
\end{description}
\end{exercise}

Next we consider the equator circle 
\[
C= \{{\bold x}\in {\Bbb S}^2: x^2 +y^2 =1\} = 
\{{\bold x}: x^2+y^2 = 1, z=0\}
\] 
in the sphere.  The stereographic projection 
$\pi(C)$ of $C$ is particularly simple---it is $C$ itself.

What happens if we first rotate the sphere and then stereographically 
project?  Say we rotate about the $y$-axis for example---by an angle $\phi$.  
If we call this rotation $R^y_\phi$, then $R^y_\phi(C) = 
\{(x\cos\phi, y, x\sin\phi): x^2+y^2 = 1\}$, and the stereographic 
projection is 
\begin{equation}\label{dirmap}
	\pi\circ R^y_\phi(C) = \left\{ \left(
				{x\cos\phi\over{1-x\sin\phi}}, 
				{y\over{1-x\sin\phi}} \right) : 
		x^2+y^2 = 1 \right\}.
\end{equation}
This set in the plane is (perhaps) not so easy to recognize.  

On page 19 of \cite{AhlCom}, however, Lars Ahlfors gave a nice way to look 
at it:  $R^y_\phi(C)$ is the intersection of a plane 
$P= \{{\bold x}: x\sin\phi - z\cos\phi = 0\}$ with ${\Bbb S}^2$.  
If we could find 
an inverse map $\pi^{-1}: {\Bbb R}^2 \to {\Bbb S}^2$, then 
\[
	\pi\circ R^y_\phi(C) = \{ {\bold a}=(a,b): 
			\pi^{-1}({\bold a}) \in P \}.
\]
Stereographic projection on ${\Bbb S}^2\backslash\{(0,0,1)\}$ does have an 
inverse (though the extended map $\bar{\pi}$ does not---why?).  The formula 
for the inverse is
\[
	\pi^{-1}({\bold a}) = {1\over{|{\bold a}|^2 + 1}}
				(2a,2b, |{\bold a}|^2 - 1).
\]
\begin{exercise} 
	Derive this formula for $\pi^{-1}$.
\end{exercise}

\noindent The statement $\pi^{-1}({\bold a}) \in P$ now translates 
into an equation: 
\[
	-2a\tan\phi\ + a^2 +b^2 = 1.
\]
Such an equation, as we know, represents a {\em circle},
\[
	(a - \tan\phi)^2 +b^2 = \sec^2\phi,
\]
with center $(\tan\phi,0)$ and radius $r = |\sec\phi|$.

If you think about what happens to $C$ when the sphere is rotated, it is 
fairly clear that the simple rotation $R^y_\phi$ is typical.  In 
particular, $\pi\circ R(C)$ should be a circle for any rotation 
$R$ of ${\Bbb S}^2$.  One way to make this precise is to decompose an 
arbitrary rotation $R$ into simple coordinate rotations like $R^y_\phi$.  
The following theorem gives such a decomposition.
\begin{theorem}\label{rthreed}
Any rotation $R$ is the composition of three rotations---one about the 
$z$-axis, one about the $y$-axis, and another about the $z$-axis.  Thus, 
there are angles $\theta$, $\phi$, and $\psi$ such that 
\[
	R= R^z_\psi\circ R^y_\phi\circ R^z_\theta.
\]
\end{theorem}

\noindent{\bf Proof.}  The simple rotations $R^z_\theta$, $R^y_\phi$, and 
$R^z_\psi$, like rotations of ${\Bbb R}^2$, can be represented by 
matrices.  For example,
\[
	R^z_\psi({\bold x}) = \left( \begin{array}{ccc} 
				\cos\psi & -\sin\psi & 0 \\
				\sin\psi & \ \ \,\cos\psi & 0 \\
				0 & 0 & 1 \end{array}\right){\bold x}
\]
for every ${\bold x}\in {\Bbb R}^3$.  Consequently, these rotations are 
{\em linear\/} transformations of ${\Bbb R}^3$.  Moreover, we show in the 
appendix that {\em all rotations are linear.\/}  Because of this, we only 
need to find $\theta$, $\phi$, and $\psi$ such that 
$R^z_\psi\circ R^y_\phi\circ R^z_\theta$ and $R$ agree on an orthonormal 
basis.

Let ${\bold e}_1 = (1,0,0)$, ${\bold e}_2 = (0,1,0)$ and 
${\bold e}_3 = (0,0,1)$.  Notice that for any $\theta$ we have 
$R^z_\theta({\bold e}_3) = {\bold e}_3$ so that 
$R^z_\psi\circ R^y_\phi\circ R^z_\theta({\bold e}_3) = 
R^z_\psi\circ R^y_\phi({\bold e}_3)$.  Thus, we want to find $\phi$ and 
$\psi$ such that, in particular, 
\begin{equation}\label{search}
R({\bold e}_3) = R^z_\psi\circ R^y_\phi({\bold e}_3) = 
	\left( \begin{array}{c} 
				-\cos\psi\sin\phi \\
				-\sin\psi\sin\phi \\
				 \cos\phi \end{array}\right)
\end{equation}
where the second equality is by direct calculation.  On the other hand, 
$R({\bold e}_3)$ is some unit vector $(u_1,u_2,u_3)$, and since 
$|u_3| \le 1$, we can find an angle $\phi$ with $\cos\phi = u_3$.  If 
$|u_3| = 1$, then $\sin\phi = 0$ and (\ref{search}) holds.  Otherwise, 
$\sin\phi \ne 0$ and $(u_1/\sin\phi)^2 + (u_2/\sin\phi)^2 = 1$.  Thus, 
for some angle $\psi$, $\cos\psi = -u_1/\sin\phi$, 
$\sin\psi = -u_2/\sin\phi$, and again (\ref{search}) holds.

We have then for any $\theta$, 
$R^z_\psi\circ R^y_\phi\circ R^z_\theta({\bold e}_3) = R({\bold e}_3)$.  
It remains to choose $\theta$ so that 
$R^z_\psi\circ R^y_\phi\circ R^z_\theta({\bold e}_j) = R({\bold e}_j)$ for 
$j = 1$ and $2$.

Let $S = R^z_\psi\circ R^y_\phi$ and note that $S^{-1}\circ R({\bold e}_3) = 
{\bold e}_3$.  Therefore, 
$S^{-1}\circ R$ induces a rotation $R_\theta$ (by some angle $\theta$) 
on ${\Bbb R}^2$.  That is, $S^{-1}\circ R({\bold e}_j) = 
R_\theta({\bold e}_j)$ for $j = 1$ and $2$.  Since $R_\theta$ 
extends to a rotation $R^z_\theta$ on ${\Bbb R}^3$ we are done.  $\Box$

\bigskip

We can now apply Theorem~\ref{rthreed} as follows.
\begin{eqnarray*}
	\pi\circ R(C) & = & \pi\circ
		R^z_\psi\circ R^y_\phi\circ R^z_\theta(C) \\
		      & = & R^z_\psi\circ \pi \circ R^y_\phi(C).
\end{eqnarray*}
This set is the circle $(a - \tan\phi)^2 +b^2 = \sec^2\phi$ rotated by 
an angle $\psi$, i.e., it is the circle with center 
$(\cos\psi \tan\phi, \sin\psi \tan\phi)$ and radius $r = |\sec\phi|$.

\begin{exercise}\label{work}
Show that any circle in ${\Bbb S}^2$ stereographically projects to a 
circle in ${\Bbb R}^2$, and any circle in ${\Bbb R}^2$ is the projection 
of a circle in ${\Bbb S}^2$.  Hint:  A circle in ${\Bbb S}^2$ is the 
intersection of a plane in ${\Bbb R}^3$ with ${\Bbb S}^2$.
\end{exercise}

We will need one other fact about the stereographic projection 
$\pi: {\Bbb S}^2 \to {\Bbb R}^2$:
\begin{lemma}\label{conformal1}
	If $C_1$ and $C_2$ are two (smooth) curves in ${\Bbb S}^2$ that 
	intersect at a point ${\bold p}$ in an angle $\gamma$, then the 
	image curves $\pi(C_1)$ and $\pi(C_2)$ intersect at $\pi({\bold p})$ 
	in the same angle $\gamma$.
\end{lemma}
This property is expressed by saying that $\pi$ is {\em conformal\/}.  
Henry Wente told us a short proof of Lemma~\ref{conformal1} which we have 
included in an appendix.  
For another proof, we refer the reader to the classic book \cite{HilGeo} by 
David Hilbert.  

\section{Symmetry}\label{SYM}

When we described the stereographic projection of a rotation of the 
equator circle ($\pi\circ R(C)$) in the last section, we ignored the 
unpleasant possibility that $C$ had been rotated onto $(0,0,1)$---where 
stereographic projection is not defined.  This happens, of course, when 
the angle $\phi$ of the preceding section is $\pi/2$.  Our analysis in 
that case is flawed since $\tan\phi$ and $\sec\phi$ are not defined, and 
in fact the image of $R(C)$ (aside from the point $(0,0,1)$) is then a 
{\em line} in ${\Bbb R}^2$.  It is one of our objectives in this section 
to address this apparent difficulty.  Our second and main objective is to 
generalize in a natural way our intuitive notion of symmetry, so that we 
can introduce the symmetry assertion of the main theorem.

The analysis of stereographic projection of circles that pass through 
$(0,0,1)$ is simple.  For example, if we take $\phi = \pi/2$, then in 
place of (\ref{dirmap}) we have
\begin{eqnarray*}
	\pi\circ R^y_{\pi/2}(C\backslash\{(1,0,0)\}) & = & \left\{ \left(
				0, 
				{y\over{1-x}} \right) : 
		x^2+y^2 = 1,\ x\ne1\right\}\\
	& = & \left\{ \left(0, 
				\pm\sqrt{1+x\over{1-x}} \right) : 
		-1\le x < 1 \right\}.	
\end{eqnarray*}
The last set is clearly the line $x=0$, since $(1+x)/(1-x)$ maps $[-1,1)$ 
monotonely onto $[0,\infty)$.  The more general cases can be handled 
similarly.  What we really wish to emphasize, however, is the following:  
{\em because circles in ${\Bbb S}^2$ that pass through $(0,0,1)$ are 
geometrically identical to other circles that do not, it is natural 
for our purposes to view straight lines in ${\Bbb R}^2$ as circles 
with (infinite radius and) one point at $\infty$.}  

For starters, this viewpoint allows Exercise~\ref{work} to make sense 
as stated.  More importantly it illustrates how sets and structures in 
${\Bbb S}^2$ can provide insight for terminology and constructions in 
${\Bbb R}^2$.  We proceed further along this line presently.

The symmetry of a circle in ${\Bbb R}^2$ is perhaps most easily described 
in terms of its center.  Given a point ${\bold p}$ in a circle $A$ with 
center ${\bold a}$, $A$ is generated by rotating ${\bold p}$ about 
${\bold a}$.  The same circle can also be generated by reflecting ${\bold p}$ 
about each of the lines through ${\bold a}$.  This latter characterization 
will be the one of interest to us.  
\begin{definition}\label{def1}
A set $A\subset {\Bbb R}^2$ has {\em Euclidean reflectional symmetry 
with respect to a point ${\bold a}\in {\Bbb R}^2$} if, for each line 
$E$ passing through ${\bold a}$, we have $\psi_E(A) = A$, where 
$\psi_E: {\Bbb R}^2\to {\Bbb R}^2$ is the reflection about $E$.  
\end{definition}
\begin{exercise}\label{euc-ex}
Show that any such set (with Euclidean reflectional symmetry) is a union 
of concentric circles with center ${\bold a}$.  
\end{exercise}

If, as we have suggested, lines should be considered simply as circles with 
infinite radius, then it is natural to ask for a definition of symmetry 
in which reflection about lines is replaced with {\em reflection about 
circles.\/}  Extrapolating directly from the definition above we might try to 
replace the family of symmetry lines passing through ${\bold a}$ with a 
family of {\em symmetry circles\/} passing through a common point 
${\bold a}$.  Unfortunately, the full geometric situation is not completely 
evident from considering the plane alone.  Again we turn to stereographic 
projection.  The inverse image of each line through ${\bold a}$ is a 
circle in ${\Bbb S}^2$ passing through $\pi^{-1}({\bold a})$ and $(0,0,1)$.  
(There are two points of intersection.)  If one then rotates slightly this 
family 
of circles in ${\Bbb S}^2$ and stereographically projects, a family of 
circles in ${\Bbb R}^2$ is obtained that pass through two distinct points 
${\bold a}_1$ and ${\bold a}_2$, as shown on the left in 
Figure~\ref{fig2}.  Note that this family also contains the 
line determined by ${\bold a}_1$ and ${\bold a}_2$.  We call the circles 
that pass through two given points ${\bold a}_1$ and ${\bold a}_2$ 
{\em Steiner symmetry circles}.  We use these circles to generalize the 
symmetry lines in Definition~\ref{def1}.

\begin{figure}[ht]
\centerline{
\begin{picture}(130,60)
	\put(88,90){${\bold p}$}
	\put(24,25){${\bold q}$}
	\put(148,115){${\bold c}$}
	\put(93,83){\circle*{3}} 
	\put(21,32){\circle*{3}} 
	\put(144,119){\circle*{3}}
	\put(0,0){\epsfxsize=2in\epsffile{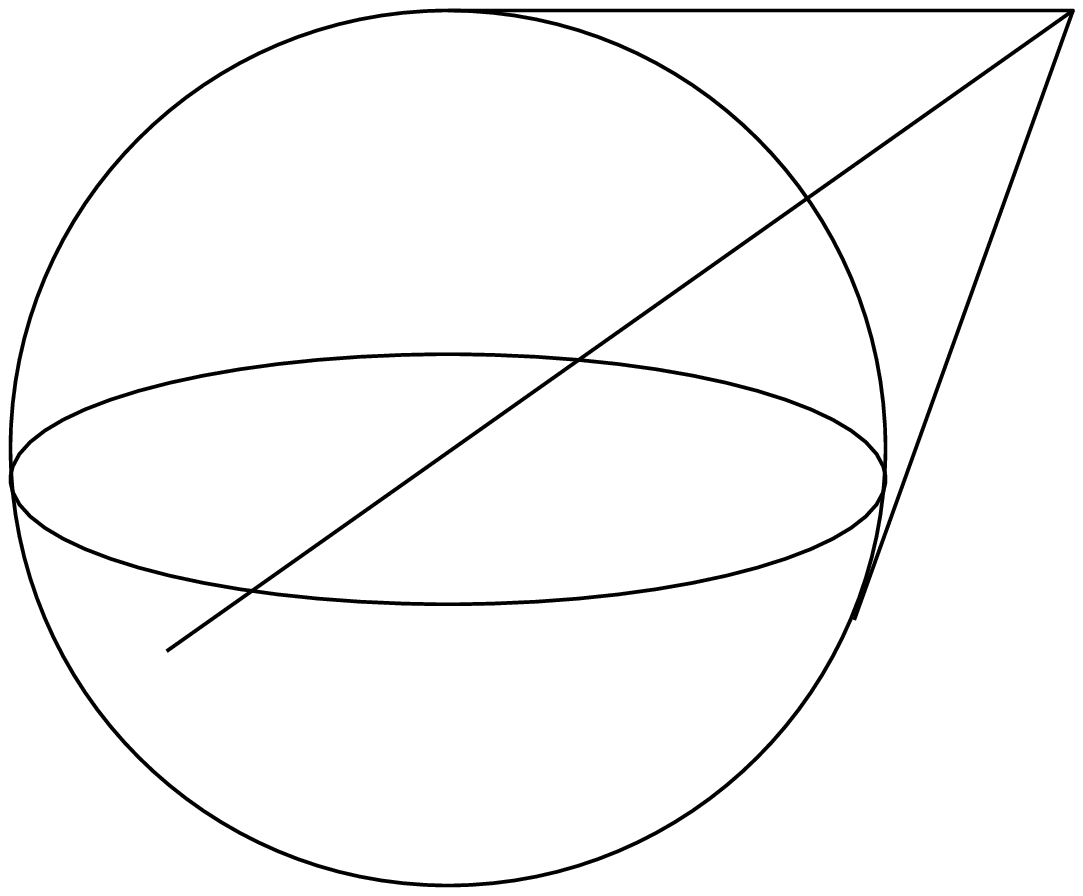}}
\end{picture}}
\centerline{ }
\caption{Reflection on ${\Bbb S}^2$.}
\label{figA}
\end{figure}

It remains to specify, for each circle $S$ passing through 
${\bold a}_1$ and ${\bold a}_2$, a transformation $\psi_S: {\Bbb R}^2 \to 
{\Bbb R}^2$ which we will call {\em reflection about $S$.\/}  Again we look 
to the sphere for intuition.  What transformation of ${\Bbb S}^2$ corresponds 
to reflection about a line in ${\Bbb R}^2$?  If $E$ is a line in 
${\Bbb R}^2\backslash\{{\bold 0}\}$, then the inverse image of $E$ is some 
circle $C$ in ${\Bbb S}^2$ passing through $(0,0,1)$.  Furthermore, there 
is a unique point ${\bold c} = (c_1,c_2,1)$ such that the segments 
connecting ${\bold c}$ to $C$ form a right circular cone (tangent to 
${\Bbb S}^2$), n.b., Figure~\ref{figA}.  The point ${\bold c}$ can now be used to define a 
transformation $\Psi_C:{\Bbb S}^2\to {\Bbb S}^2$:  
\begin{quote}
	For each point ${\bold p}\in{\Bbb S}^2$, the line determined 
	by ${\bold p}$ and ${\bold c}$ intersects ${\Bbb S}^2$ in a 
	set $\{{\bold p},{\bold q}\}$.  We set $\Psi_C({\bold p})={\bold q}$.  
\end{quote} 
Note that this definition is much 
like the geometric definition of stereographic projection.  
The reader can check (and we will show below) that $\Psi_C$ corresponds to 
the reflection $\psi_E$ in the sense that 
$\psi_E = \pi\circ\Psi_C\circ\pi^{-1}$.  This construction also works if 
$E$ passes through ${\bold 0}$---though in that case the transformation 
$\Psi_C$ is simply given by reflection about the plane determined by $C$.

Notice that the fact $C$ passes through $(0,0,1)$ is not required for the 
geometric definitions of $\Psi_C$ to make sense.  That is, for any circle 
$C$ in ${\Bbb S}^2$, if $C$ is not a great circle, it defines a cone 
point\footnote{Or horizon point.  See \cite{HilGeo} for other interesting 
properties of this point.  In particular, $\bar{\pi}({\bold c})$ is the center 
of $\pi(C)$.} ${\bold c}$, and the definition above gives a transformation 
$\Psi_C$ of ${\Bbb S}^2$ that is geometrically identical (modulo rotation) 
to one that corresponds to Euclidean reflection.  For great circles we 
use the alternative construction.  

Going back to ${\Bbb R}^2$, we may start with any circle (or straight 
line) $S$, take $C=\pi^{-1}(S)$ and apply the construction described 
above to obtain a transformation $\psi_S = \pi\circ\Psi_C\circ\pi^{-1}$ 
of ${\Bbb R}^2\backslash\{ {\bold a} = \bar{\pi}({\bold c})\}$.  
This is the transformation we call {\em reflection about 
$S$.}

\begin{figure}[ht]
\centerline{
\begin{picture}(300,160)
	\put(0,0){\epsffile{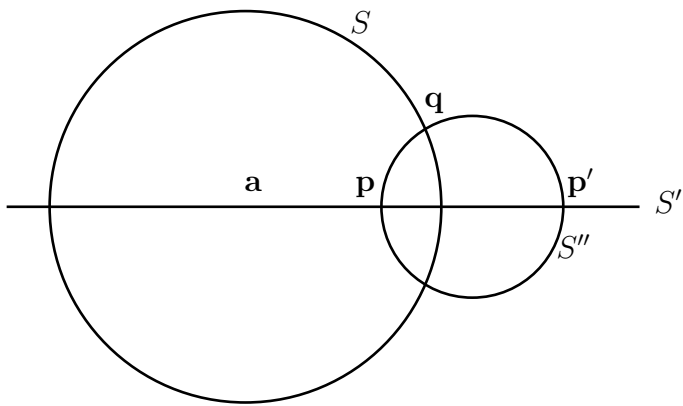}}
	\put(130,140){$S$}
	\put(90,80){${\bold a}$}
	\put(208,55){$S''$}
	\put(245,72){$S'$}
	\put(132,80){${\bold p}$}
	\put(212,80){${\bold p'}$}
	\put(158,112){${\bold q}$}
\end{picture}}
\centerline{ }
\caption{Mapping Circles}
\label{fig1}
\end{figure}

\begin{figure}[ht]
\centerline{
\begin{picture}(350,150)
	\put(15,84){${\bold a}_1$}
	\put(63,42){${\bold a}_2$}
	\put(35,82){\circle*{3}}
	\put(62,55){\circle*{3}}
	\put(180,100){$l$}
	\put(0,0){\epsfxsize=5in\epsffile{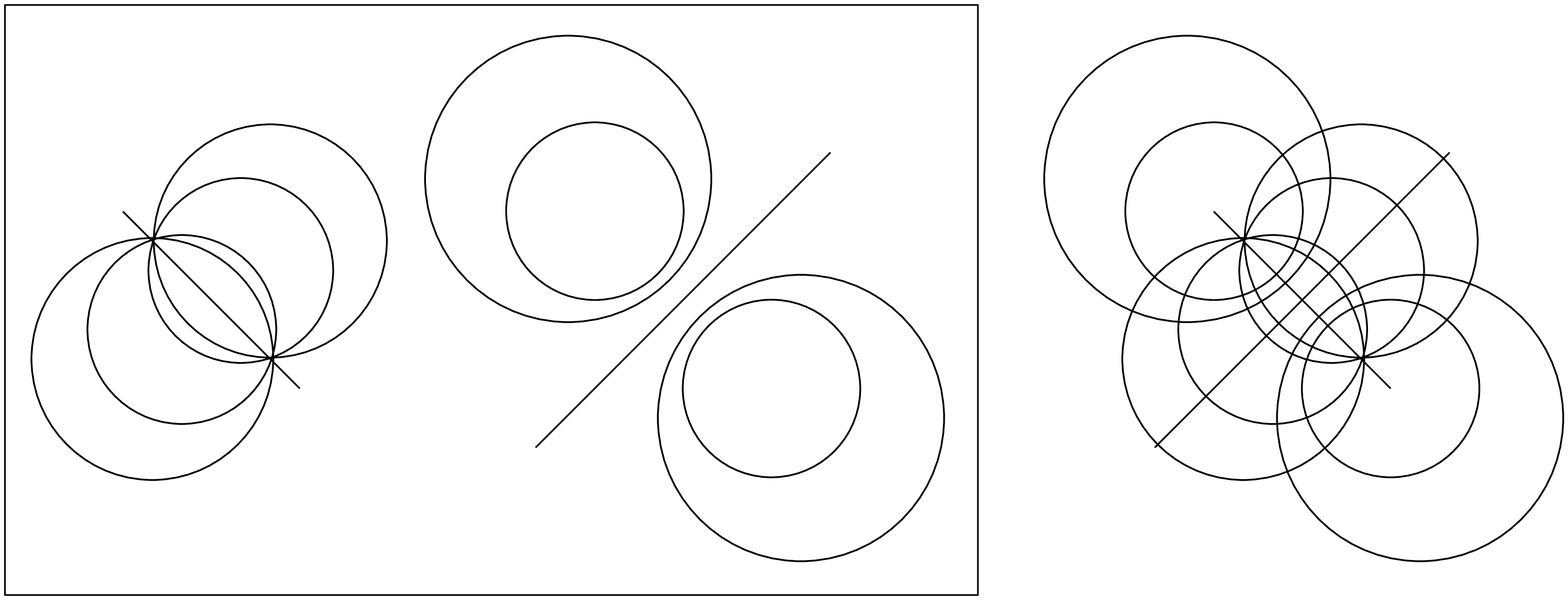}}
\end{picture}}
\centerline{ }
\caption{Steiner circles and circles of Apollonius.}
\label{fig2}
\end{figure}

\begin{figure}[ht]
\centerline{
\begin{picture}(440,120)
	\put(40,0){$Q(0)$}
	\put(200,0){$Q(\pi/8)$}
	\put(345,0){$Q(\pi/4)$}
	\put(0,0){\epsfxsize=2in\epsffile{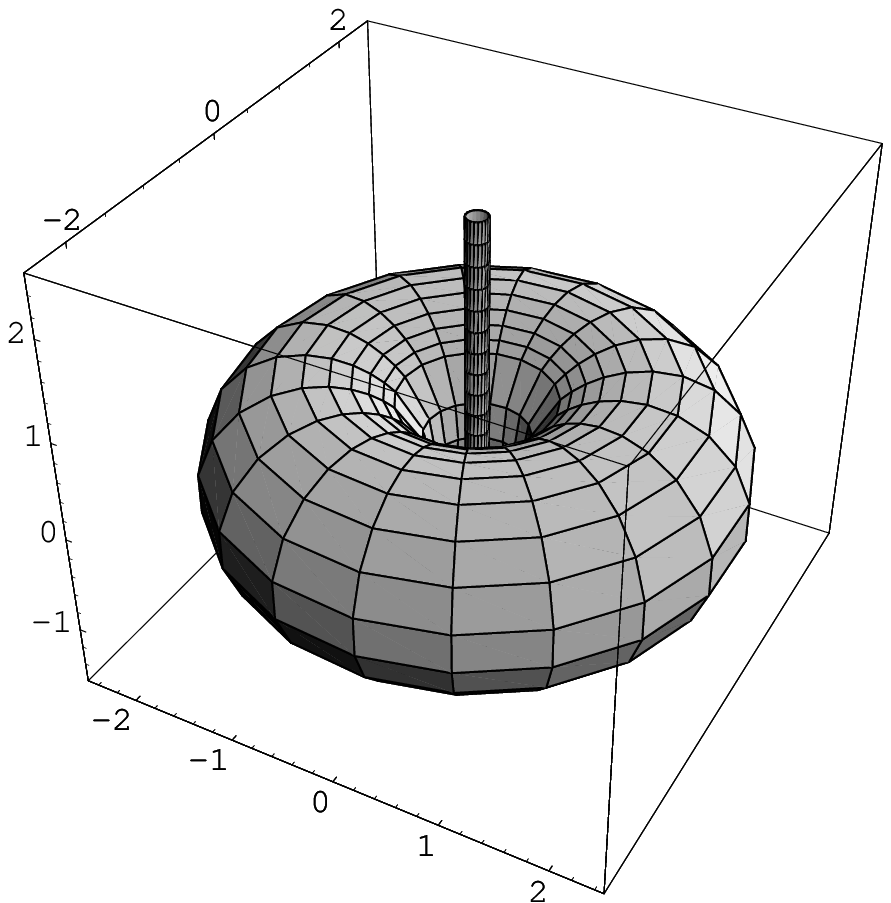}}
	\put(150,0){\epsfxsize=2.3in\epsffile{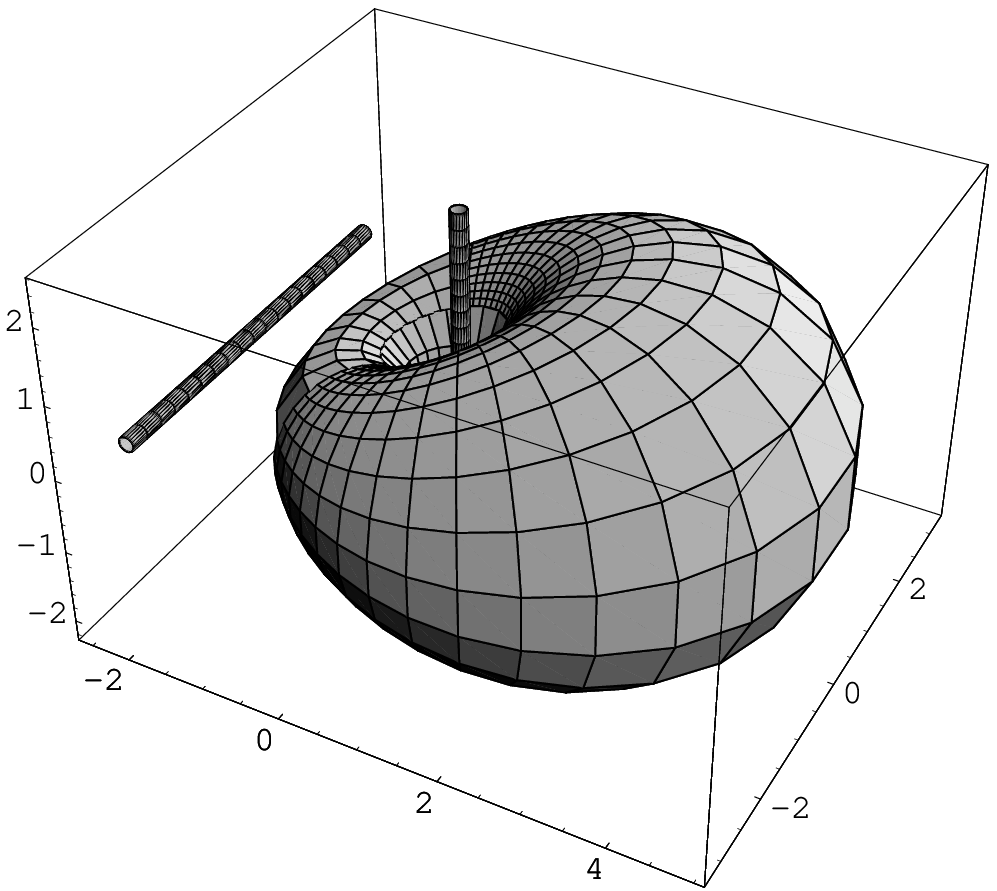}}
	\put(335,15){\epsfxsize=1in\epsffile{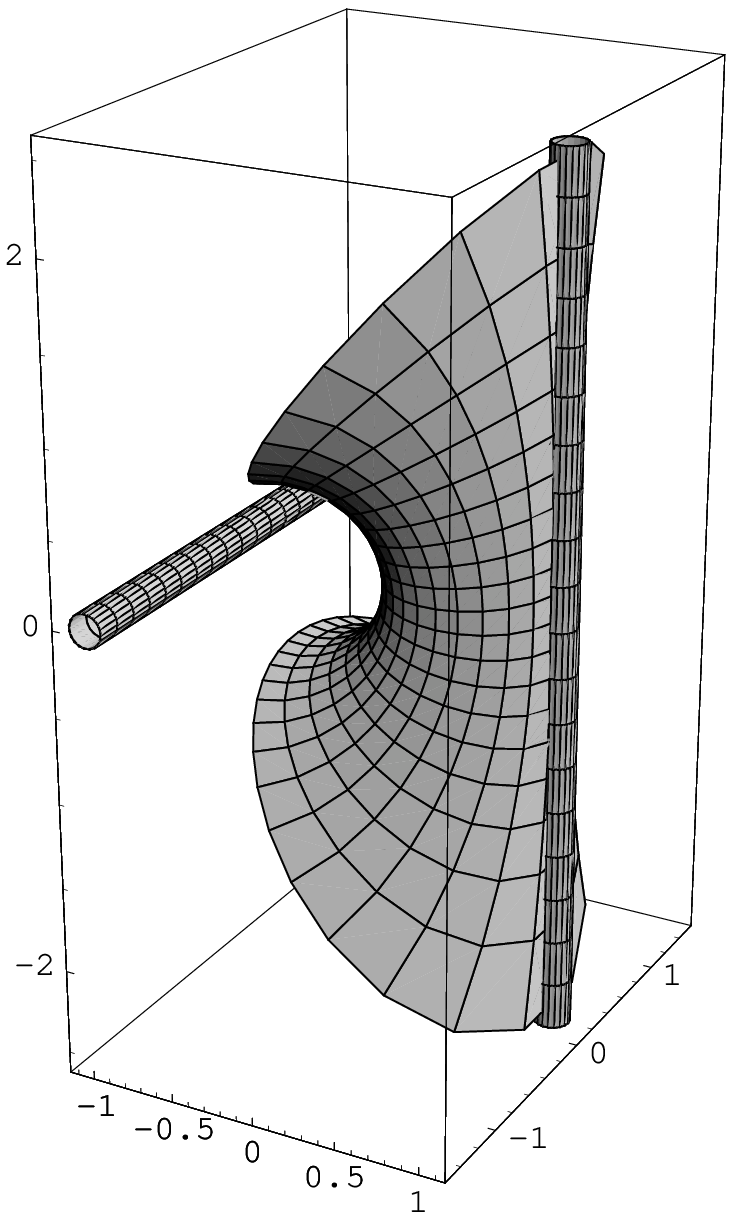}}
\end{picture}}
\centerline{ }
\caption{Stereographic projections and lines of centers.}
\label{figcliff}
\end{figure}

We proceed to derive a formula for 
$\psi_S: {\Bbb R}^2\backslash\{ {\bold a} \} \to {\Bbb R}^2$.  
Notice first that if $C'$ is a circle in ${\Bbb S}^2$ that meets $C$ 
in right angles at points ${\bold p}_1$ and ${\bold p}_2$, then $C'$ 
lies in the plane determined by ${\bold c}$, ${\bold p}_1$ and ${\bold p}_2$, 
and one sees from this that $\Psi_C(C') = C'$.  It follows moreover, since 
$\pi$ is an angle and circle preserving transformation, that any circle 
$S'$ in ${\Bbb R}^2$ which is orthogonal to $S$ is mapped by $\psi_S$ 
into itself.  Thus, consider a point ${\bold p}$ inside $S$ as shown in 
Figure~\ref{fig1}.  Let 
$S'$ be the line determined by the center of $S$ and ${\bold p}$.  There 
is also a circle $S''$ which passes through ${\bold p}$, meets $S$ 
orthogonally, and has its center on $S'$.  Since 
$\psi_S(S'\backslash\{{\bold a}\}) = S'\backslash\{{\bold a}\}$, 
$\psi_S(S'') = S''$, and $\psi_S({\bold p}) \ne {\bold p}$, it follows 
that $\psi_S({\bold p}) = {\bold p'}$ is the other point 
of intersection of $S'$ and $S''$.  
Moreover, if ${\bold a}$ is the center of $S$ and ${\bold q}\in S\cap S''$, 
then triangles ${\bold a}{\bold q}{\bold p}$ and 
${\bold a}{\bold q}{\bold p'}$ are similar.  It follows that 
$|{\bold p}-{\bold a}||{\bold p'}-{\bold a}| = \rho^2$ where $\rho$ is 
the radius of $S$.  The same reasoning applies if ${\bold p}$ lies outside 
of $S$, and we obtain the formula
\begin{equation}\label{formula}
	{\bold p'} = \psi_S({\bold p}) = 
		\rho^2 {{\bold p}-{\bold a}\over{|{\bold p}-{\bold a}|^2} }
			+ {\bold a}
\end{equation}
for reflection about the circle $S$ in ${\Bbb R}^2$ with center ${\bold a}$ 
and radius $\rho$.  This same discussion applied to the case when $C$ 
passes through $(0,0,1)$ and $S=E$ is a straight line provides a proof 
of the assertion made above that $\Psi_C$ corresponds to $\psi_E$.

\thispagestyle{empty}

\medskip

Finally we have the following

\medskip

\noindent{\bf Definition}\ {\em  
A set $A\subset {\Bbb R}^2$ has {\em generalized 
reflectional symmetry} if there are two distinct points 
${\bold a}_1$ and ${\bold a}_2$ such that for each Steiner circle $S$ passing 
through ${\bold a}_1$ and ${\bold a}_2$, $\psi_S(A) = A$ where 
$\psi_S$ 
is reflection about $S$.  }

\bigskip

We may allow one of the points ${\bold a}_1$ or ${\bold a}_2$ to 
be at $\infty$,  in which case the circles $S$ are all the lines 
passing through the other point.  Furthermore, it can be shown (in 
analogy to Exercise~\ref{euc-ex}) that {\em if $A$ has generalized 
reflectional symmetry, then $A$ is a union of circles, each of which 
is orthogonal to all the circles through ${\bold a}_1$ and ${\bold a}_2$.}  
We will not need this fact, but it is a special case of an assertion 
proved in \cite{McCDis}.  
These circles are called {\em circles of Apollonius\/}, and 
once it is known that they are circles, it is easy to show the following.  
(See Figure~\ref{fig2}.)
\begin{lemma}The circles of Apollonius determined by the Steiner circles 
through ${\bold a}_1$ and ${\bold a}_2$ in ${\Bbb R}^2$ are disjoint and 
are in one to one 
correspondence with their centers which comprise the line $E$ passing through 
${\bold a}_1$ and ${\bold a}_2$ except for the segment between 
${\bold a}_1$ and ${\bold a}_2$.  Let ${\bold m}_0 = 
({\bold a}_1+{\bold a}_2)/2$.  The circle of Apollonius with center 
${\bold a}\in {\Bbb R}^2$ has radius $r = \sqrt{d^2 - \rho_0^2}$ where 
$d= |{\bold a} - {\bold m}_0|$ and 
$\rho_0 = |{\bold a}_1 - {\bold a}_2|/2$.  
\end{lemma}
The centers of the Steiner circles also form a line $l$, and once we know 
the point ${\bold m}_0 = ({\bold a}_1+{\bold a}_2)/2$ on $l$ and the 
reference distance $\rho_0 = |{\bold a}_1 - {\bold a}_2|/2$, we can express 
generalized reflectional symmetry without reference to 
${\bold a}_1$ or ${\bold a}_2$ as follows.
\begin{definition}\label{def-ref}
A set $Q$ has {\em generalized reflectional symmetry along a line $l$} if, 
for each point ${\bold a}\in l$, we have $\psi_S(Q) = Q$, where 
$\psi_S$ is given 
by {\rm (\ref{formula})} with $\rho = \sqrt{d^2 +\rho_0^2}$ and 
$d = |{\bold a} - {\bold m}_0|$.
\end{definition}

Notice finally that formula (\ref{formula}) makes sense for points 
${\bold a}$ and ${\bold p}$ in ${\Bbb R}^3$ and gives a generalization 
of reflection about circles to reflection about spheres.  Thus, this 
 statement of generalized symmetry can be applied to sets $Q\subset
{\Bbb R}^3$.  

In the next section we use this formulation and take $Q$ to be a 
stereographic projection of the Clifford torus.

\section{Stereographic projection of the Clifford Torus}\label{MR}
The unit sphere in ${\Bbb R}^4$ is the three-dimensional space 
\[
	{\Bbb S}^3 = \{ {\bold x}=(x,y,z,w): |{\bold x}| = 1\}.
\]
Because we are used to visualizing things that are described by three 
Euclidean coordinates (i.e., things in ${\Bbb R}^3$), it is often 
difficult to see what objects look like in ${\Bbb S}^3$.  For this reason, 
it is convenient to use a stereographic projection $\pi: {\Bbb S}^3\backslash 
\{(0,0,0,1)\} \to {\Bbb R}^3$.  The formula for such a map is similar 
to the one for ${\Bbb S}^2$:
\[
	\pi({\bold x}) = {1\over{1-w}}(x,y,z),
\]
and a similar geometric description applies as well.

We are interested in a particular geometric object in ${\Bbb S}^3$ called the 
{\em Clifford torus}:
\[
	{\cal C} = \{ {\bold x}: x^2 +y^2 = 1/2 = z^2 + w^2 \}.
\]
The stereographic projection $\pi(\cliff)$ of the Clifford torus is 
particularly nice because of its {\em symmetry}.
\begin{exercise}\label{sym.ex}
\begin{description}
\item[(i)] Show that $\pi(\cliff)$ is rotationally symmetric with respect 
	to the $z$-axis in ${\Bbb R}^3 = \{ (x,y,z,0) \}$.
\item[(ii)] What is the intersection of $\pi(\cliff)$ with the 
	$x,z$-plane?
\end{description}
\end{exercise}
From Exercise~\ref{sym.ex} it is clear that $\pi(\cliff)$ is described by 
its intersection with the half planes $\Pi_\theta = 
\{ (r\cos\theta, r\sin\theta, z): r > 0 \}$.  In fact, it is enough 
to know only $\pi(\cliff) \cap \Pi_0$.  

As with circles in ${\Bbb S}^2$, rotating the three-sphere (i.e., moving 
$\cliff$ around in ${\Bbb S}^3$) changes the stereographic 
projection.  
Since the rotated surface is geometrically identical to $\cliff$ however, one 
might expect that some kind of symmetry of the projection is preserved.  
In fact, we show the following.
\begin{theorem}\label{mainthm}
Let $R$ 
be any rotation of ${\Bbb S}^3 \subset {\Bbb R}^4$.  
The surface $Q = \pi\circ R(\cliff)$ has generalized symmetry in the sense 
described in the last section.
\end{theorem}

We will, as we did in \S1, consider first a particular rotation and 
then show that that rotation is typical via a decomposition formula 
(Theorem~\ref{decomp2}) for general rotations.  With hindsight from the 
decomposition formula, we consider the rotation $R^{xw}_\psi$ of the 
$x,w$-plane corresponding to the matrix
\[
	\left(\begin{array}{cccc}
		\cos\psi & 0 & 0 & -\sin\psi  \\
		0        & 1 & 0 & 0          \\
		0        & 0 & 1 & 0          \\
		\sin\psi & 0 & 0 & \cos\psi    \end{array}\right).
\]
\begin{proposition}  \label{prop}
Let $Q = Q(\psi) = \pi \circ R^{xw}_\psi({\cal C})$. 
\begin{description}
	\item[(i)] If $\psi = 0$, then $Q = Q_0$ has generalized reflectional 
symmetry along the $z$-axis.  
	\item[(ii)] If $0<|\psi| \le \pi/2$, then $Q$ has generalized 
reflectional symmetry along the vertical line $(x,y) = (\tan\psi,0)$ and 
along the horizontal line $(x,z) = (-\cot\psi,0)$.
	\item[(iii)]  Let $R$ be the rotation of ${\Bbb R}^3$ about the 
$x$-axis by $\pi/2$.  Then
	\begin{description}
		\item[(a)] $Q(\psi+\pi/2) = R(Q(\psi))$, and
		\item[(b)] $Q(\psi+\pi) = Q(\psi)$, for all $\psi$.
	\end{description}
\end{description}
In particular, $Q=Q(\psi)$ has generalized reflectional symmetry 
for all $\psi$.
\end{proposition}
\begin{theorem}\label{decomp2}
Any rotation $R$ of ${\Bbb S}^3\subset{\Bbb R}^4$ is a composition 
\begin{equation}\label{form2}
	R = R_0 \circ R^{xw}_\psi \circ R^{zw}_\phi \circ R^{xy}_\theta
\end{equation}
where $R_0$ is the trivial extension\footnote{\label{fnum}For a precise 
definition see the appendix on rotations.} to ${\Bbb R}^4$ of a rotation of 
${\Bbb R}^3 = \{ (x,y,z,0) \}$ and 
$R^{xy}_\theta$ is the trivial extension$^{\ref{fnum}}$ to ${\Bbb R}^4$ 
of a rotation of ${\Bbb R}^2 = \{ (x,y,0,0) \}$ etc..
\end{theorem}

Before we prove Proposition~\ref{prop} and Theorem~\ref{decomp2} we will 
show that they imply Theorem~\ref{mainthm}.  It is easy to check (see 
Lemma~\ref{switches} below) that 
rotations of the $x,y$ and $z,w$-planes leave ${\cal C}$ invariant.  
Furthermore, if $S$ is any set in ${\Bbb S}^3 \backslash\{ (0,0,0,1) \}$ 
and $R_0$ is a rotation of $\{ (x,y,z,0) \}$ as above, then 
$\pi\circ R_0(S) = R_0\circ\pi(S)$.  To see this, let $R_0(x,y,z,0) = 
(x',y',z',0)$ and note that 
\begin{eqnarray*}
\pi\circ R_0(x,y,z,w) & = & \pi(x',y',z',w) \\
	& = & {1\over{1-w}}(x',y',z',0) \\
	& = & R_0\left({1\over{1-w}} (x,y,z,0)\right) \\
	& = & R_0\circ\pi(x,y,z,w).
\end{eqnarray*}

Thus, using the decomposition formula (\ref{form2}) we have
\begin{eqnarray*}
	\pi\circ R({\cal C}) & = & \pi\circ R_0 \circ 
		R^{xw}_\psi \circ R^{zw}_\phi \circ R^{xy}_\theta ({\cal C}) \\
			& = & R_0 \circ \pi \circ R^{xw}_\psi ({\cal C}).
\end{eqnarray*}
The set $\pi \circ R^{xw}_\psi ({\cal C})$ is described in 
Proposition~\ref{prop}.  The rotation $R_0$ only changes that description by 
a rigid rotation in ${\Bbb R}^3$.  In particular, $\pi\circ R({\cal C})$ has 
generalized reflectional symmetry.  $\Box$

\bigskip

\thispagestyle{empty}

\noindent{\bf Proof of Proposition~\ref{prop}.}  Part (i) follows from 
Exercise~\ref{sym.ex} where one finds that Definition~\ref{def-ref} is 
satisfied with ${\bold m}_0 = {\bold 0} \in {\Bbb R}^3$ and $\rho_0 = 1$.  
Therefore, the sphere of symmetry with center ${\bold a} = (0,0,c)$ which 
we denote by $S = S_\rho({\bold a})$ has radius 
\begin{equation}\label{rho}
	\rho = \sqrt{1+c^2}.
\end{equation}

For $0< |\psi| \le \pi/2$ we will show that $C = \pi^{-1}(S)$ determines 
a transformation $\Psi_C$ of ${\Bbb S}^3$.  This construction is analogous 
to the discussion in \S\ref{SYM} of circles $C \subset {\Bbb S}^2$.  Moreover, 
we will again have the correspondence $\psi_S = \pi\circ\Psi_C\circ\pi^{-1}$ 
where $\psi_S$ is the reflection about $S$.  Furthermore, a geometrically 
identical transformation $\Psi_{\tilde{C}}$ will be determined by 
$\tilde{C} = R^{xw}_\psi(C)$, and $\Psi_{\tilde{C}}$ will correspond to 
reflection about the {\em sphere\/} $\tilde{S} = \pi(\tilde{C})$.  The 
spheres $\tilde{S}$ thus corresponding to the spheres $S= S_\rho({\bold a})$ 
will be symmetry spheres for $Q = Q(\psi)$ that satisfy 
Definition~\ref{def-ref}.  Of course, at this point we do not even know that 
$\tilde{S}$ is a sphere.  We give now a precise higher dimensional version 
of Exercise~\ref{work} which will establish this fact.  For the statement 
we use the notation $\underline{{\bold n}} = (n_1,n_2,n_3)$ when 
${\bold n} = (n_1,n_2,n_3,n_4)$ and the notation $({\bold n},n_4) = 
(n_1,n_2,n_3,n_4)$ when ${\bold n} = (n_1,n_2,n_3)$.
\begin{lemma}\label{sph2}
Let $\Pi = \{ {\bold x}=(x,y,z,w): {\bold n}\cdot {\bold x} = e\}$ be a 
three-plane in ${\Bbb R}^4$.  If $(0,0,0,1)\notin \Pi$, then 
\begin{equation}\label{3p2s}
\pi(\Pi\cap {\Bbb S}^3) = \left\{ {\bold a} = (a,b,c): \left|{\bold a} - 
		{\underline{{\bold n}}\over{n_4 - e}}\right|^2 = 
		{n_4+e\over{n_4-e}} + 
		{|\underline{{\bold n}}|^2\over{(n_4 - e)^2}} \right\}
\end{equation}
(which is a sphere).  If $(0,0,0,1)\in \Pi$, then 
\begin{equation}\label{3p2p}
\pi(\Pi\cap {\Bbb S}^3\backslash\{(0,0,0,1)\}) = 
		\{ {\bold a}: \underline{{\bold n}}\cdot {\bold a} = n_4 \}
\end{equation}
(which is a plane).

On the other hand, let 
$S = \{{\bold a}: |{\bold a} - {\bold a}_0|^2 = \rho^2 \}$ be a sphere in 
${\Bbb R}^3$.  Then $\pi^{-1}(S) = \Pi \cap {\Bbb S}^3$ where 
\begin{equation}\label{s23p}
\Pi = \{ {\bold x}: (-2{\bold a}_0, \rho^2 - |{\bold a}_0|^2 + 1)\cdot 
		{\bold x} = \rho^2 - |{\bold a}_0|^2 - 1\}.
\end{equation}
Let $P = \{{\bold a}: {\bold n}\cdot {\bold a} = e \}$ be a plane in 
${\Bbb R}^3$.  Then $\pi^{-1}(P) \cup\{(0,0,0,1)\} = \Pi \cap {\Bbb S}^3$ 
where 
\begin{equation}\label{p23p}
\Pi = \{ {\bold x}: ({\bold n},e)\cdot {\bold x} = e\}.
\end{equation}
\end{lemma}
Since we have given explicit equations, Lemma~\ref{sph2} follows from 
simple substitution using the formulas for $\pi$ and $\pi^{-1}$, and we 
omit the proof.  Note that equations (\ref{3p2s}) and (\ref{3p2p}) allow 
degenerate cases corresponding to $\Pi\cap {\Bbb S}^3 \subset \{(0,0,0,1)\}$.  
In our applications below however, we will know that $\Pi\cap {\Bbb S}^3$ is 
nontrivial.

\medskip

Recall that $S = S_\rho({\bold a})$ is a sphere of symmetry for $Q_0$.  
One sees from (\ref{s23p}) and (\ref{rho}) that $C = \pi^{-1}(S) = 
\Pi \cap {\Bbb S}^3$ where $\Pi= \{ {\bold x}: ({\bold a}, -1)\cdot {\bold x} 
= 0 \}$.  Note that $\Pi$ passes through ${\bold 0}\in {\Bbb R}^4$, i.e., 
$C$ is a {\em great sphere} in ${\Bbb S}^3$.  Let ${\bold n}= ({\bold a},-1)$ 
be the normal to $\Pi$.  We consider the reflection $\Psi_C$ of ${\Bbb S}^3$ 
about $\Pi$ defined by
\begin{equation}\label{4ref}
\Psi_C({\bold x}) = {\bold x} - {2{\bold x}\cdot {\bold n} \over{
			|{\bold n}|^2}} {\bold n}.
\end{equation}
To see that $\psi_S = \pi\circ\Psi_C\circ\pi^{-1}$, we extend the discussion 
of Figure~\ref{fig1} in \S\ref{SYM}.  Let $S' = \{ {\bold a}': {\bold n}' 
\cdot ({\bold a}' - {\bold a}) = 0 \}$ be a plane orthogonal to $S$.  
According to (\ref{p23p}) we have $\pi^{-1}(S') \cup \{(0,0,0,1)\} = 
\Pi' \cap {\Bbb S}^3$ where $\Pi' = \{ {\bold x}: ({\bold n}',{\bold n}'\cdot 
{\bold a}) \cdot {\bold x} = {\bold n}'\cdot {\bold a} \}$.  Since 
$({\bold n}',{\bold n}'\cdot {\bold a}) \cdot {\bold n} = 
({\bold n}',{\bold n}'\cdot {\bold a}) \cdot ({\bold a},-1) = 0$, we see 
that $\Psi_C(\Pi') = \Pi'$.  Consequently, 
$\pi\circ\Psi_C\circ\pi^{-1} (S'\backslash\{{\bold a}\}) = 
S'\backslash\{{\bold a}\}$.  It follows similarly that 
$\pi\circ\Psi_C\circ\pi^{-1} (S'') = S''$ for any sphere $S''$ orthogonal to 
$S$.

It then follows that formula (\ref{formula}) gives the value of 
$\pi\circ\Psi_C\circ\pi^{-1} ({\bold p})$, and hence that $\psi_S = 
\pi\circ\Psi_C\circ\pi^{-1}$.  To see this, we can apply the discussion of 
Figure~\ref{fig1} in \S\ref{SYM} where we interpret $S$, $S'$ and $S''$ as 
$S_\rho({\bold a})$, a plane (orthogonal to the paper), and a sphere 
respectively.  Technically, we should also introduce the plane of the 
paper $S'''$ which can be used to show that ${\bold p}'\in S'''$.

There is nothing special about $S = S_\rho({\bold a})$ in this reasoning, 
except that its inverse image is a great sphere.  We have actually shown 
the following.
\begin{lemma}\label{corresp}
Let $S$ be any sphere or plane in ${\Bbb R}^3$ whose inverse image 
$C = \pi^{-1}(S)$ is determined by a three-plane $\Pi = \{ {\bold x}\cdot 
{\bold n} = 0 \}$ through ${\bold 0}\in {\Bbb R}^4$.  Then the reflection 
$\psi_S$ about $S$ is given by $\psi_S = \pi\circ\Psi_C\circ\pi^{-1}$ where 
$\Psi_C$ is the reflection about $\Pi$ given by {\rm (\ref{4ref})}.
\end{lemma}

We are now in a position to finish the proof of Proposition~\ref{prop}.  
The rotation $\tilde{C} = R^{xw}_\psi(C)$ is also a great sphere in 
${\Bbb S}^3$, and for $0<|\psi|\le \pi/2$, the plane $\tilde{\Pi} = 
R^{xw}_\psi(\Pi) = 
\{ {\bold x}: (\sin\psi,0,c,-\cos\psi)\cdot {\bold x} = 0 \}$ does not 
contain $(0,0,0,1)$.  Thus, we have from (\ref{3p2s}) that $\pi(\tilde{C})$ 
is the sphere
\[
\tilde{S} = \{ \tilde{{\bold a}}: 
	|\tilde{{\bold a}} - (\tan\psi,0,c/\cos\psi)|^2 = 
	(1+c^2)/\cos^2\psi \}.
\]
If we set $\tilde{{\bold m}}_0 = (\tan\psi,0,0)$ and $\tilde{\rho_0} = 
1/\cos^2\psi$, we see that $Q = \pi\circ R^{xw}_\psi({\cal C})$ satisfies 
Definition~\ref{def-ref} since
\begin{eqnarray*}
	\psi_{\tilde{S}}(Q) & = & \pi\circ\Psi_{\tilde{C}}\circ\pi^{-1}\ \circ
				\ \pi\circ R^{xw}_\psi({\cal C}) \\
		& = & \pi\circ R^{xw}_\psi\circ \Psi_C \circ R^{xw}_{-\psi} 
			\ \circ \ R^{xw}_\psi({\cal C}) \\
		& = & \pi\circ R^{xw}_\psi\ \circ\ \pi^{-1}\circ\pi \circ 
			\Psi_C \circ \pi^{-1}\circ\pi({\cal C}) \\
		& = & \pi\circ R^{xw}_\psi\ \circ\ \pi^{-1}\circ\psi_S 
			\circ\pi({\cal C}) \\
		& = & \pi\circ R^{xw}_\psi({\cal C}) = Q.
\end{eqnarray*}
We have established that $Q$ has a vertical line $(x,y) = (\tan\psi,0)$ 
of generalized reflectional symmetry for $0 < |\psi| \le \pi/2$.

\medskip

The planes $P_\phi = 
\{ {\bold a}: (-\sin\phi,\cos\phi,0) \cdot {\bold a} = 0\}$ 
are also planes of reflective symmetry for $Q_0$.  Applying 
Lemma~\ref{corresp} and Lemma~\ref{sph2} much as we have done above with 
the spheres $S_\rho({\bold a})$ we find that for $0 < |\psi| \le \pi/2$ 
the spheres 
\begin{eqnarray*}
	\tilde{S} & = & \pi\circ R^{xw}_\psi \circ\pi^{-1}(P_\phi) \\
		  & = & \{ {\bold a}: |{\bold a} - 
		(-\cot\psi,\cot\phi/\sin\psi,0)|^2 = 
		(1+\cot^2\phi)/\sin^2\psi \} 
\end{eqnarray*}
are spheres of symmetry along the horizontal line $(x,z) = (-\cot\psi,0)$ 
which satisfy Definition~\ref{def-ref} with $\tilde{{\bold m}}_0 = 
(-\cot\psi,0,0)$ and $\tilde{\rho}_0 = 1/\sin^2\psi$.  This finishes the 
proof of statement (ii).

The first identity in statement (iii) follows from explicit calculation 
and the following observation.  
\begin{lemma}\label{switches} 
For any function $f:{\Bbb R}^4 \to {\Bbb R}^4$ we 
have $\{ f(x,y,z,w): {\bold x} \in {\cal C}\} = 
\{ f(x',y',z',w'): {\bold x} \in {\cal C}\}$ where 
$x'$ is $\pm x$ and $y'$ is $\pm y$ (or possibly $x'$ is $\pm y$ and 
$y'$ is $\pm x$), and  similarly $z'$ is $\pm z$ and $w'$ is $\pm w$ 
(or possibly $z'$ is $\pm w$ and $w'$ is $\pm z$).  
\end{lemma}

Statement (iiib) follows from (iiia). $\Box$

\bigskip

\noindent{\bf Proof of Theorem~\ref{decomp2}.}  We use again the fact 
that rotations are precisely those linear transformations that 
correspond to orthogonal matrices of determinant 1.  Let $M$ be the 
matrix representing $R$.  Let $N=N(\theta,\phi,\psi)$ be the unknown 
matrix representing $R^{xw}_\psi \circ R^{zw}_\phi \circ R^{xy}_\theta$, 
and let $N_0$ be the unknown matrix representing $R_0$.  We then need to 
show $M = N_0N$.  

We know that $N_0$ has the form
\begin{equation}\label{matform}
N_0 = \left(\begin{array}{cccc}
			&              &   & 0 \\
			& {\huge\star} &   & 0 \\
			&              &   & 0 \\
		      0 & 0            & 0 & 1 \end{array}\right)
\end{equation}
where $\Large\star$ represents a rotation matrix for ${\Bbb R}^3$.  From 
this we see that the last row of $N_0N$ and the last row of $N$ are the 
same.  This last row is given by $N^{\rm T}{\bold e}_4$ (where 
${\rm T}$ indicates 
the transpose), and we need to have $M^{\rm T}{\bold e}_4 = 
N^{\rm T}{\bold e}_4$.  
\begin{lemma}
There exist angles $\theta$, $\phi$, and $\psi$ such that 
$M^{\rm T}{\bold e}_4 = N^{\rm T}{\bold e}_4$. 
\end{lemma}
\noindent{\bf Proof.}  Since $M$ is orthogonal, so is $M^{-1} = M^{\rm T}$.  
Therefore, the columns of $M^{\rm T}$ (i.e., the rows of $M$) form an 
orthonormal basis---see \cite[pp. 127--129]{CurLin}.  In particular, 
$M^{\rm T}{\bold e}_4 = (m_{41},m_{42},m_{43},m_{44})$ is some unit 
vector.

On the other hand, by direct calculation we see that 
\[
N^{\rm T}{\bold e}_4 = 
(\sin\psi\cos\theta,-\sin\psi\sin\theta,\cos\psi\sin\phi,\cos\psi\cos\phi).
\]  
Since $m_{41}^2 + m_{42}^2 \le 1$, there is some angle $\psi$ with 
$\sin^2\psi = m_{41}^2 + m_{42}^2$.  If $\sin\psi \ne 0$, then we can 
find $\theta$ with $\cos\theta = m_{41}/\sin\psi$ and 
$\sin\theta = - m_{42}/\sin\psi$.  If $\sin\psi = 0$, then 
$m_{41} = m_{42} = 0$.  In either case, our choice of $\theta$ and 
$\psi$ implies that the first two coordinates of $M^{\rm T}{\bold e}_4$ 
and $N^{\rm T}{\bold e}_4$ agree.  Since 
$m_{43}^2 + m_{44}^2 = 1- m_{41}^2 + m_{42}^2 = \cos^2\theta$, we can 
choose $\phi$, much as we chose $\theta$, and have the last two 
coordinates match.  $\Box$

\medskip

To prove Theorem~\ref{decomp2} it remains to specify $N_0$.  Let 
${\bold m}_j = (m_{j1}, m_{j2}, m_{j3}, m_{j4})$ be the $j$th row of 
$M$ and ${\bold n}_j$ be the $j$th row of $N$ for $j=1,2,3,4$, and 
consider 
\begin{eqnarray*}
MN^{-1} & = & \left(\begin{array}{cccc}
			{\bold m}_1 \\
			{\bold m}_2 \\
			{\bold m}_3 \\
			{\bold n}_4  \end{array}\right) 
		\left({\bold n}_1,{\bold n}_2,{\bold n}_3,{\bold m}_4\right) \\
	& = & \left(\begin{array}{cccc}
			&              &   & {\bold m}_1 \cdot {\bold m}_4 \\
			& {\huge\star} &   & {\bold m}_2 \cdot {\bold m}_4 \\
			&              &   & {\bold m}_3 \cdot {\bold m}_4 \\
	{\bold n}_4 \cdot {\bold n}_1 & {\bold n}_4 \cdot {\bold n}_2  
	& {\bold n}_4 \cdot {\bold n}_3 & {\bold n}_4 \cdot {\bold m}_4 
		\end{array}\right) \\
	& = & \left(\begin{array}{cccc}
			&              &   & 0 \\
			& {\huge\star} &   & 0 \\
			&              &   & 0 \\
		      0 & 0            & 0 &  {\bold m}_4 \cdot {\bold m}_4 
		\end{array}\right)
\end{eqnarray*}
which is of the form (\ref{matform}).  Thus, we let $N_0$ be the rotation 
matrix on the right and clearly we have $M = N_0N$.  $\Box$

\section*{Epilogue}  We observed in \S2 that stereographic projections 
of rotations of the equator circle in ${\Bbb S}^2$ are circles.  
Exercise~\ref{work} points out that other circles in ${\Bbb S}^2$ have 
this property as well.  That is, we have {\em not\/} characterized the 
rotations of the equator circle $C$.  Nevertheless, $C$ and its rotations 
(the great circles) are ``balanced'' on the surface of ${\Bbb S}^2$ in a 
way that the other circles are not.  This ``balance'' is expressed 
precisely by saying the {\em geodesic curvature\/} is zero or simply that 
these curves are {\em geodesics\/}.  It is this balance that justifies the 
specific attention we have given to the great circles.

In a similar way, the Clifford torus is ``balanced'' in the three-sphere 
${\Bbb S}^3$, because its {\em mean curvature} is zero, i.e., it is a 
{\em minimal surface\/}.  While it can be shown that every geodesic curve 
in ${\Bbb S}^2$ is (part of) a great circle, there is a great variety of 
minimal surfaces in ${\Bbb S}^3$.  There are spheres and tori and surfaces 
of genus two (two holed tori), etc..  In fact, Lawson \cite{LawCom} has 
given examples of closed minimal surfaces in ${\Bbb S}^3$ of every 
{\em topological genus\/}, i.e., tori with any number of holes.  
On the other hand, Bryant \cite{BryDua} has shown that every embedded 
(non self-intersecting) minimal surface in ${\Bbb S}^3$ that is 
topologically spherical stereographically projects to a standard 
Euclidean sphere.  A similar characterization for minimal tori is not 
known, but it is believed that, up to a rigid rotation of ${\Bbb S}^3$, 
the Clifford torus ${\cal C}$ is the unique embedded minimal torus.

For this reason, symmetry properties of ${\cal C}$ as shown above are 
of great interest.  We remark finally that the only closed surfaces in 
${\Bbb R}^3$ possessing the symmetry shown above for stereographic 
projections of ${\cal C}$ are topologically spherical or toroidal 
\cite{McCDis}.  This symmetry, moreover, has geometric consequences 
as well, and we hope to give, in another paper, an elementary 
introduction to the curvature of surfaces in ${\Bbb S}^3$ and prove 
that {\em ${\cal C}$ is the unique minimal torus possessing such 
symmetry.}

\appendix

\section*{Appendix: \ Rotations}
We seek below to give an intuitive introduction to the family of 
{\em rotations\/} of Euclidean space ${\Bbb R}^n$, $n\ge 3$.  
Our starting point is with the distance preserving transformations 
$T:{\Bbb R}^n \to {\Bbb R}^n$ which satisfy
\begin{equation}\label{rm}
|T({\bold x}) - T({\bold y})| = |{\bold x} - {\bold y}|.
\end{equation}
We refer to all such transformations as {\em rigid motions}, and our 
objective is to determine which rigid motions should be called 
{\em rotations}.

Notice first of all that {\em translations\/} are rigid motions.  That is, 
for any fixed vector ${\bold a}$, the transformation defined by 
$T({\bold x}) = 
{\bold x} + {\bold a}$ satisfies (\ref{rm}).  
It is a fundamental algebraic fact that up to a 
translation every rigid motion is {\em linear}.
\begin{theorem}\label{rigid.lem}
If $T_0:{\Bbb R}^n \to {\Bbb R}^n$ is a rigid motion, then $T$ defined by 
$T({\bold x}) = T_0({\bold x}) - T_0({\bold 0})$ is a linear transformation.
\end{theorem}

\noindent{\bf Proof.\footnote{This discussion is considered in a more 
general setting in \cite[pp. 1--6]{AusInt}.}}  
Recall that $T$ is 
linear if $T(a{\bold x}) = aT({\bold x})$ and 
$T({\bold x}) + T({\bold y}) = T({\bold x} + {\bold y})$ for all 
${\bold x},{\bold y}\in {\Bbb R}^n$ and $a\in {\Bbb R}^1$.  Recall also 
the {\em triangle inequality\/} in ${\Bbb R}^n$:  
\begin{quote}{\em 
	$|{\bold x} + {\bold y}| \le |{\bold x}| + |{\bold y}|$ 
	with equality only if ${\bold x} = \lambda{\bold y}$ for some 
	$\lambda \ge 0$. }
\end{quote}
See \cite[Exercise 1-2]{SpiCal} for a proof.  

Because $T$ preserves distance and fixes the origin,
\begin{equation}\label{norms}
	|T({\bold x})| = |{\bold x}|,\ 
	|T(a{\bold x})| = |a|\,|{\bold x}|,\ 
	|T(a{\bold x}) - T({\bold x})| = |a-1|\,|{\bold x}|. 
\end{equation}
On the other hand, by the triangle inequality 
\begin{equation}\label{ti}
	|T({\bold x}) - T(a{\bold x})| + |T(a{\bold x})| \ge |T({\bold x})|.
\end{equation}
If $0\le a\le 1$, one can use (\ref{norms}) to check that equality holds 
in (\ref{ti}).  Consequently, for some $\lambda \ge 0$ 
\[
	T({\bold x}) - T(a{\bold x}) = \lambda T(a{\bold x}),
\]
or $T({\bold x}) = (1+\lambda) T(a{\bold x})$.  Taking the norm of both 
sides we see that $T(a{\bold x}) = a T({\bold x})$.  

By exchanging $|T(a{\bold x})|$ and $|T({\bold x})|$ in (\ref{ti}), and 
following the same line of reasoning, one sees that 
$T(a{\bold x}) = a T({\bold x})$ also for $1<a$.

Finally, if $a<0$, the same reasoning applied to the inequality 
\[
	|T(a{\bold x})| + |T({\bold x})| \ge |T({\bold x}) - T(a{\bold x})|
\]
yields again that $T(a{\bold x}) = a T({\bold x})$.

\medskip

Next consider $T({\bold x}+{\bold y})$.  In fact, let $a<0$ and 
note that
\[
	|T({\bold x}+a{\bold y}) - T({\bold x})| + 
		|T({\bold x}) - T({\bold x}+{\bold y})| \ge 
		|T({\bold x}+a{\bold y}) - T({\bold x}+{\bold y})|.
\]
The left side is $|a{\bold y}| + |{\bold y}| = (1-a){\bold y}$, and the 
right side is $|a{\bold y} - {\bold y}| = (1-a){\bold y}$.  Since they 
are equal, there exists $\lambda \ge 0$ with 
$T({\bold x}+a{\bold y}) - T({\bold x}) = 
\lambda (T({\bold x}) - T({\bold x}+{\bold y}))$.  It is easy to see 
that $\lambda = -a$, so 
\[
	T({\bold x}+a{\bold y}) = T({\bold x}) + 
			a(T({\bold x}+{\bold y}) - T({\bold x})).
\]
Subtracting $aT({\bold y}) = T(a{\bold y})$ from both sides and 
rearranging we get 
\[ 
	a[T({\bold x}+{\bold y}) -(T({\bold x})+ T({\bold y}))] = 
		T({\bold x}+a{\bold y}) - T(a{\bold y}) - T({\bold x}).
\]
Therefore, 
\begin{eqnarray*}
	|a|\,|T({\bold x}+{\bold y}) -(T({\bold x})+ T({\bold y}))| & \le &
		|T({\bold x}+a{\bold y}) - T(a{\bold y})| + |T({\bold x})| \\
			& = & 2|{\bold x}|.
\end{eqnarray*}
Notice that the right side is a fixed value, but $|a|$ on the left may 
be taken as {\em large}\/ as we like.  The only way the inequality can 
continue to hold is if 
\[
	|T({\bold x}+{\bold y}) -(T({\bold x})+ T({\bold y}))| = 0,
\]
i.e., $T({\bold x}+{\bold y}) = T({\bold x})+ T({\bold y})$.  $\Box$

\bigskip

From now on, we assume our rigid motions satisfy 
\begin{equation}\label{normalize}
T({\bold 0}) = {\bold 0}.
\end{equation}
Thus, we are only going to consider rotations about the origin---since 
other rotations only differ from these by a translation.  Furthermore, we 
can refer to the matrix $M$ which corresponds to a rigid motion $T$ (which 
will be the matrix of $T$ with respect to the standard basis unless stated 
otherwise).

From Theorem~\ref{rigid.lem} we can easily prove 
\begin{corollary}\label{basis-preserved}
If ${\bold e}_1,\ldots, {\bold e}_n$ form an orthonormal basis, then 
$T({\bold e}_1),\ldots, T({\bold e}_n)$ form an 
orthonormal basis as well.
\end{corollary}

\noindent{\bf Proof.}  Once linearity 
is established, the preservation of orthonormal bases 
follows by expressing the inner product in terms of the norm 
$|{\bold x}| = \sqrt{{\bold x}\cdot {\bold x}}$.  In fact, 
\[
	{\bold x}\cdot {\bold y} = {1\over{2}} ( |{\bold x} + {\bold y}|^2 
		- |{\bold x}|^2 - |{\bold y}|^2). 
\]
It follows, using the linearity, that 
$T({\bold x}) \cdot T({\bold y}) = {\bold x}\cdot {\bold y}$.  
In particular, $T({\bold e}_i)\cdot T({\bold e}_j) = 
{\bold e}_i \cdot {\bold e}_j 
= \delta_{ij}$ where ${\bold e}_1,\ldots,{\bold e}_n$ form 
an orthonormal basis.  $\Box$

\medskip

\noindent Any linear transformation that preserves orthonormality of 
bases is called an {\em orthogonal} transformation.  The main properties 
of orthogonal transformations and 
the matrices that represent them may be found in 
\cite[pp. 127--129]{CurLin}.  In particular, the inverse matrix $M^{-1}$ 
that represents $T^{-1}$ is the transpose matrix $M^T$ of $M$.  It follows 
from this and the product formula for determinants that $({\rm det\,}M)^2 
= 1$, or since $M$ is a real matrix, that ${\rm det\,}M = \pm 1$.  As 
pointed out in the introduction, the additional condition
\begin{equation}\label{det1}
{\rm det\,}M = 1
\end{equation}
is often used to distinguish $T$ as a rotation.

\medskip

So far we have used the intuitive condition (\ref{rm}) and the normalization 
(\ref{normalize}) to derive some algebraic facts.  We now return to our 
intuition concerning rotations and ask for a precise condition which, in 
conjunction with (\ref{rm}), will express our intuitive idea of what 
defines a rotation.  This is not so easy, but we should keep in mind that 
such a condition is likely to be equivalent to (\ref{det1}).

\begin{figure}[hb]
\centerline{
\begin{picture}(240,120)
	\thicklines
	\put(12,36){\circle*{8}}
	\put(0,37){${\bold 0}$}
	\put(43,75){${\bold e}_2$}
	\put(62,25){${\bold e}_1$}
	\put(12,0){\epsfxsize=3in\epsffile{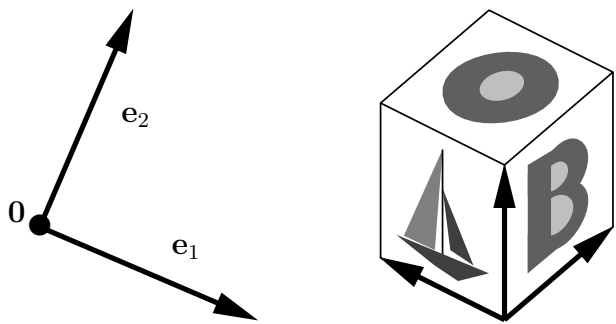}}
\end{picture}}
\centerline{ }
\caption{Moving orthonormal bases.}
\label{fig3}
\end{figure}

Our first approach might be to give an intuitive (yet precise) definition 
of {\em orientation}, and then try to connect a condition concerning 
{\em orientation preserving} transformations with (\ref{det1}).  The 
reader may be surprised to find, as we were, that our physical intuition 
concerning this approach is limited to two dimensions.  To see this, take 
a piece of paper and draw an orthonormal basis on it 
dark enough so that you can see it from the back side of the paper (see 
Figure~\ref{fig3}).  
Now draw a dot on the desk to represent the origin ${\bold 0}$.  Any 
rigid motion of the plane (i.e., the paper) which fixes the origin 
must map ${\bold e}_1$ to 
some unit vector $(\cos\theta,\sin\theta)$---which you can draw on the 
desk.  Now there are two obvious ways to rigidly move the paper so that 
it lies flat on the desk and ${\bold e}_1$ lies on top of 
$(\cos\theta,\sin\theta)$.  Intuitively, if the paper is facing up, the motion 
is a rotation.  If the paper is facing down, it is not.  That is, whether 
or not a rigid motion is a rotation is determined by how we ``orient'' 
the paper before placing it on the desk.  

It is very difficult however (if not impossible) for us to rigidly move 
a physical representation of ${\Bbb R}^3$ (like a wooden block) so as 
to change its orientation.  

It turns out that the only way to define {\em orientation} of bases for 
${\Bbb R}^3$ is, one way or another, to append an additional dimension.  
One feels, however, that we {\em do} have an intuitive idea of what 
constitutes a rotation of ${\Bbb R}^3$ independent of additional dimensions.

A second approach might be based on the idea of a {\em rotation axis.}  
Indeed, every rotation of ${\Bbb R}^3$ has $1$ as an eigenvalue so that it 
does have a fixed vector, ${\bold x} = R({\bold x})$, which can be used to 
define an axis of rotation (See \cite[pg. 291, Corollary 33.3]{CurLin}).  
Unfortunately, there is no such rotation axis for non-trivial rotations of 
${\Bbb R}^2$, and there need not be one for rotations of ${\Bbb R}^4$.  

A third approach (since we are getting frustrated) could be to use the 
idea of {\em decomposition} as in Theorems~\ref{rthreed} and 
\ref{decomp2}---except in reverse.  That is, we could define an 
{\em elementary rotation} to be a rotation of just one coordinate two-plane, 
i.e., a transformation $R^{kj}_\psi$ corresponding to a matrix of the 
form
\[
\begin{array}{c} k\to \\  \\ j\to 
			\\ \end{array}
\left(\begin{array}{ccccc} I_{k-1} & & & & \\
			& \cos\psi & & -\sin\psi & \\
			& & I_{j-k-1} & & \\
			& \sin\psi & & \ \cos\psi & \\
			& & & & I_{n-j} \end{array}\right)
\]
where $I_m$ denotes an $m\times m$ identity matrix and there are zeros filling 
all the spaces.  Then we could define a {\em rotation} to be a composition 
of elementary rotations.  After pondering this, however, it is not at all 
clear that a composition of rotations 
should be a 
rotation.  In fact, the decomposition in Theorem~\ref{rthreed} is not 
intuitively a rotation of ${\Bbb R}^3$---it is the composition of three 
rotations, one executed after another {\em in time.}  And here is the 
key.  A {\em rotation} is a transformation which can be realized as a 
physical rigid motion (parameterized by time) that is the {\em same} motion at 
each instant of time.  To make this statement precise is fairly easy.
\begin{definition}\label{rots}
A rigid motion $R$ is a {\em rotation} if there is a smoothly parameterized 
family of rigid motions $R_0(t)$ such that $R_0(0) = {\rm id}_{{\Bbb R}^n}$ 
and, 
for each $m= 1,2,3,\ldots$,
\[
R_0(1/m)^m = \underbrace{R_0(1/m)\circ \cdots 
\circ R_0(1/m)}_{m\ {\rm times}}
= R.
\]
\end{definition}
Notice that $R_0(t)$ for $t\in(0,1)$ is not explicitly required to be a 
rotation (so the definition is not circular).  One might be worried 
however that the definition allows transformations of determinant $-1$ 
which we don't want as rotations in ${\Bbb R}^2$.  It is easy to show 
that this does not happen, but it turns out that the most difficult thing 
to see is that condition (\ref{det1}) does not allow transformations that 
Definition~\ref{rots} excludes.  Nevertheless, we have the following.
\begin{proposition}\label{prop2}
If $R:{\Bbb R}^n \to {\Bbb R}^n$ is a rigid motion represented by the 
matrix $M$, then the following are equivalent.
\begin{enumerate}
\item[{\rm (i)}] $R$ is a rotation.
\item[{\rm (ii)}] There exist elementary rotations $R_1,R_2,\ldots,R_k$ such that 
$R= R_1\circ\cdots\circ R_k$.
\item[{\rm (iii)}] ${\rm det\,} M = 1$.  
\item[{\rm (iv)}] With respect to some basis $R$ is represented by a matrix of 
the form
\[ \left(\begin{array}{cccc} I_k & & & \\
		& R_{\theta_1} & & \\
		& & \ddots & \\
		& & & R_{\theta_{(n-k)/2}} \end{array}\right)
\]
where $I_k$ is a $k\times k$ identity matrix and 
$R_{\theta_1},\ldots,R_{\theta_{(n-k)/2}}$ 
are $2\times 2$ rotation matrices.
\end{enumerate}
\end{proposition}

Notice that condition (ii) was mentioned above as a condition that came 
somewhat short in expressing our intuitive idea of a rotation.  It is 
however a useful condition, of which Theorems~\ref{rthreed} and \ref{decomp2} 
are particular instances, so we have included it.  We will use the following 
exercise and lemma to show that it is part of the equivalence.
\begin{exercise}\label{verify.ex}
Show that elementary rotations are rotations and have determinant $1$.
\end{exercise}
\begin{lemma}\label{al1}
Given any two vectors ${\bold v},{\bold w}\in{\Bbb R}^n$ of the same 
length, there is 
a rotation $Q$, which is a composition of elementary rotations 
$R^{kl}_\theta$, such that $Q{\bold v} = {\bold w}$.
\end{lemma}

\noindent{\bf Proof.}  We first note that it is enough to prove the 
lemma for ${\bold v} = {\bold e}_n$ and ${\bold w} = {\bold u}$ an 
arbitrary {\em unit} vector.  To see this simply note that 
$\tilde{\bold v} = {\bold v}/|{\bold v}|$ and 
$\tilde{\bold w} = {\bold w}/|{\bold w}|$ are unit vectors.  Thus, 
if we can find $Q_1$ and $Q_2$ (compositions of elementary rotations) 
with $Q_1({\bold e}_n) = \tilde{\bold v}$ and 
$Q_2({\bold e}_n) = \tilde{\bold w}$, we can take $Q= Q_2\circ Q_1^{-1}$ 
and it is easily checked that $Q({\bold v}) = {\bold w}$.

We prove that ${\bold e}_n$ can be ``coordinate rotated'' to ${\bold u}$ 
by induction.  The initial case, $n=2$ follows from the fact that in 
${\Bbb R}^2$ any unit vector ${\bold u}$ can be represented by 
$(\cos\theta,\sin\theta)$ for some angle $\theta$.  

For $n>2$, let ${\bold u}= (u_1,\ldots,u_n)$ and $u_n = \cos\phi$.  It 
follows that for some ${\bold v}= (v_1,\ldots,v_{n-1}) \in {\Bbb R}^{n-1}$
\[
	R_\phi^{1n}({\bold e}_n) = (v_1,\ldots,v_{n-1}, u_n).
\]
Also, by Exercise~\ref{verify.ex}, $R_\phi^{1n}$ preserves length, so 
$|{\bold v}| = |{\bold w}|$ where ${\bold w}= (u_1,\ldots,u_{n-1})$.  By 
induction, there is a composition $Q$ of elementary rotations of 
${\Bbb R}^{n-1}$ such that $Q({\bold v}) = {\bold w}$.  Notice 
that $Q$ extends to a composition of elementary rotations of ${\Bbb R}^n$, 
and we have $Q\circ R_\phi^{1n}({\bold e}_n) = {\bold u}$.  This 
completes the induction and the proof of Lemma~\ref{al1}. $\Box$

\medskip

\noindent{\bf Proof of Proposition~\ref{prop2}.}  That (i) implies (iii) 
follows from the product formula for determinants applied to $R = 
R_0(1/2) \circ R_0(1/2)$.  That (iii) implies (iv) is Theorem~30.5 
in \cite[pg. 270]{CurLin} where rotations are viewed as orthogonal 
transformations of determinant $1$.  Condition (i) follows from (iv) by 
taking $R_0(t)$ to be the transformation corresponding (in the same basis) 
to the matrix 
\[ \left(\begin{array}{cccc} I_k & & & \\
		& R_{t\theta_1} & & \\
		& & \ddots & \\
		& & & R_{t\theta_{(n-k)/2}} \end{array}\right).
\]

Condition (ii) we deal with separately.  It is clear from the product 
formula for determinants and Exercise~\ref{verify.ex} that (ii) implies 
(iii).  We obtain the reverse implication by induction.  As discussed 
above, if $n=2$, then $R({\bold e}_1)$ must be some unit vector 
$(\cos\theta,\sin\theta)$.  By orthogonality 
$R(e_2) = \pm(-\sin\theta,\cos\theta)$.  Only the $+$ sign is 
compatible with (iii).  

For $n\ge 3$, let $M$ be the orthogonal matrix representing $R$.  
Let ${\bold v} = R({\bold e}_n)$.  According to 
Lemma~\ref{al1}, there is a composition $Q$ of elementary rotations 
such that $Q{\bold v} = {\bold e}_n$.  Thus, $Q\circ R$, fixes 
${\bold e}_n$, and if $N$ is the matrix representing $Q$ we 
have 
\begin{equation}\label{last}
	NM = \left(\begin{array}{cccc} 
			&   &   &  0  \\
			& U &   &  \vdots \\
			&   &   &  0  \\
		      0 & \cdots & 0 & 1  \end{array}\right)
\end{equation}
where $U$ is an $(n-1)\times(n-1)$ rotation matrix.  By 
Exercise~\ref{verify.ex} and the multiplication 
formula for determinants ${\rm det\,}U = {\rm det\,}(NM) = 
{\rm det\,} M = 1$.  Therefore, by induction $U = U_1\cdots U_k$ for some 
elementary rotation matrices $U_1,\ldots,U_k$.  These rotations extend 
linearly to elementary rotations of ${\Bbb R}^n$ represented by 
matrices $N_j$ of the form (\ref{last}) with $U_j$ in place of $U$.
Hence, $M= N^{-1}N_1\cdots N_k$ is a product of elementary rotation 
matrices.  $\Box$

\medskip

We note that the construction just used to extend the elementary rotation 
matrices can be used to extend any rotation to a rotation on a higher 
dimensional space.  To be precise, if $R$ is a rotation of ${\Bbb R}^k$ 
with standard basis elements ${\bold e}_1,\ldots,{\bold e}_k$, $k<n$, and 
$J = \{j_1<\cdots <j_k\}$ is a subset of $k$ indices from $\{1,\ldots,n\}$, 
then one obtains a rotation $\tilde{R}$ of ${\Bbb R}^n$ with standard basis 
elements $\tilde{{\bold e}}_1,\ldots,\tilde{{\bold e}}_n$ by setting 
$\tilde{R}(\tilde{{\bold e}}_j) =\tilde{{\bold e}}_j$ if $j\notin J$ and 
$\tilde{R}(\tilde{{\bold e}}_{j_l}) = \sum b_m \tilde{{\bold e}}_{j_m}$ where 
$R({\bold e}_{l}) = \sum b_m {\bold e}_m$.  The rotation $\tilde{R}$ 
is called a {\em trivial extension} of $R$.

\section*{Appendix:\ Conformality of Stereographic Projection}

Here we give a short proof of Lemma~\ref{conformal1}.  Let us first 
assume that the two curves $C_1$ and $C_2$ intersect at the point 
${\bold p} = (0,0,-1)$ in an angle $\gamma(C_1,C_2,{\bold p})$.  Let 
$T_j$ be a unit tangent vector to $C_j$ at $p_j$ for $j=1$ and $2$, and   
let $\tilde{C}_j$ be the intersection (circle) of the plane $\Pi_j$ 
containing $T_j$ and $(0,0,1)$ with ${\Bbb S}^2$.  Clearly we have 
$\gamma(C_1,C_2,{\bold p})= \gamma(\tilde{C}_1,\tilde{C}_2,{\bold p})$ 
is the angle between $\Pi_1$ and $\Pi_2$.

On the other hand, $\pi({\bold p}) = (0,0)$, and $\pi(\tilde{C}_j)$ is the 
intersection (line) of $\Pi_j$ with the $x,y$-plane for $j = 1$ and $2$.  
Thus, the angle of intersection of $\pi(\tilde{C}_1)$ and 
$\pi(\tilde{C}_2)$ is again the angle between $\Pi_1$ and $\Pi_2$.  We have 
shown that $\pi$ is conformal at ${\bold p} = (0,0,-1)$.

If the point of intersection ${\bold p}$ is any point other than 
$(0,0,-1)$ in ${\Bbb S}^2\backslash\{(0,0,1)\}$, we again take tangents 
$T_j$ to $C_j$ at ${\bold p}$ and let $\tilde{C}_j$ be the intersection 
(circle) determined by the plane $\Pi_j$ containing $T_j$ and $(0,0,-1)$ 
for $j=1$ and $2$.  Here $\pi(\tilde{C}_1)$ and $\pi(\tilde{C}_2)$ are 
circles that intersect in two points ${\bold q}_1 = \pi({\bold p})$ and 
${\bold q}_2 = (0,0)$.  Since circles (on the sphere and in the plane) 
intersect in equal angles at their two points of intersection, we have 
\begin{eqnarray*}
\gamma(\pi(C_1),\pi(C_2),\pi({\bold p})) & = & \gamma(\pi(\tilde{C}_1),\pi(\tilde{C}_2),\pi({\bold p})) \\ 
	& = & \gamma(\pi(\tilde{C}_1),\pi(\tilde{C}_2),(0,0)) \\
	& = & \gamma(\tilde{C}_1,\tilde{C}_2,(0,0,-1)) \\ 
	& = & \gamma(\tilde{C}_1,\tilde{C}_2,{\bold p}) \\
	& = & \gamma(C_1,C_2,{\bold p}).  \Box
\end{eqnarray*}

\hskip-\parindent
John McCuan\\Mathematics Department\\University of California,
Berkeley\\Berkeley, CA 94720\\johnm@math.berkeley.edu 


\begin{thebibliography}{10}

\bibitem{AhlCom}
Lars~V. Ahlfors.
\newblock {\em Complex Analysis}.
\newblock McGraw-Hill, New York, 1979.

\bibitem{AusInt}
Louis Auslander and Robert MacKenzie.
\newblock {\em Introduction to Differentiable Manifolds}.
\newblock McGraw-Hill, New York, 1963.

\bibitem{BryDua}
R.~Bryant.
\newblock A duality theorem for {W}illmore surfaces.
\newblock {\em J. Differential Geometry}, 20:23--53, 1984.

\bibitem{CurLin}
Charles~W. Curtis.
\newblock {\em Linear Algebra: An Introductory Approach}.
\newblock Springer-Verlag, New York, 1984.

\bibitem{HilGeo}
D.~Hilbert and S.~Cohn-Vossen.
\newblock {\em Geometry and the Imagination}.
\newblock Chelsea, New York, 1956.

\bibitem{HsuMin}
Lucas Hsu, Rob Kusner, and John Sullivan.
\newblock Minimizing the squared mean curvature integral for surfaces in space
  forms.
\newblock {\em Experimental Mathematics}, 1(3):191--207, 1992.

\bibitem{LawCom}
H.B. Lawson.
\newblock Complete minimal surfaces in {${\Bbb S}^3$}.
\newblock {\em Ann. Math.}, 92:335--374, 1970.

\bibitem{McCDis}
John McCuan.
\newblock {\em Symmetry via Spherical Reflection and Spanning Drops in a
  Wedge}.
\newblock PhD thesis, Stanford, 1995.

\bibitem{SpiCal}
Michael Spivak.
\newblock {\em Calculus on Manifolds}.
\newblock Addison-Wesley, New York, 1965.

\bibitem{WilRie}
T.J. Willmore.
\newblock {\em Riemannian Geometry}.
\newblock Oxford, New York, 1993.

\bibitem{YauSem}
Shing~Tung Yau.
\newblock {\em Seminar on Differential Geometry}.
\newblock Number 102 in Annals of Mathematics Studies. Princeton University
  Press, Princeton, N.J., 1982.

\bibitem{ZulCha}
Louis Zulli.
\newblock Charting the 3-sphere---an exposition for undergraduates.
\newblock {\em Am. Math. Monthly}, 103(3):221--229, 1996.

\end{thebibliography}
\end{document}